\documentclass[twoside,10pt,a4paper,leqno]{article}
\usepackage{latexsym}
\usepackage{amsmath,amsfonts,amscd,amssymb,bbm}
\usepackage[all]{xypic}
\usepackage[latin1]{inputenc}
\parskip=\medskipamount
\newtheorem{Satz}{Satz}%[section]
\newtheorem{Thm}[Satz]{Theorem}
\newtheorem{Ex}[Satz]{Example}
\newtheorem{Def}[Satz]{Definition}

\newtheorem{Lemma}[Satz]{Lemma}
\newtheorem{Prop}[Satz]{Proposition}

\newtheorem{Rem}[Satz]{Remark}
\newtheorem{Conj}[Satz]{Conjecture}
\def\itboxx#1{\ifvmode\indent\fi\makebox[2em][r]{\rmn(#1)} }
\newcommand{\rmn}{\fontshape{n}\fontseries{m}\selectfont\rm}
\newcommand{\rk}{\operatorname{rk}}

\newcommand{\To}{\longrightarrow}

\newcommand{\ds}{\displaystyle}

\newcommand{\id}{\operatorname{id}}

\newcommand{\ok}{\hfill $\Box$}

\begin{document}
\ \\
\vspace*{5mm}
\begin{center}{\large \bf
Group actions of prime order on local normal rings
}\\[5mm]
{\it Franz Kir\'{a}ly
\footnote{Franz.Kiraly@uni-ulm.de}
 and Werner Lütkebohmert}
\footnote{Werner.Luetkebohmert@uni-ulm.de}
\\[2mm]
Dept. of Pure Mathematics, University of Ulm, 89069 Ulm, Germany\\[2mm]
14 April 2011
\end{center}
\vspace*{5mm}

In the theory of singularities, an important class of singularities is built by the famous Hirzebruch-Jung singularities. They arise by dividing out a finite cyclic group action on a smooth surface. The resolution of these singularities is well understood and have nice arithmetic properties related to continued fractions; cf.  \cite{Hirzebruch} and \cite{Jung}.\\
 One can also look at such group actions from a purely algebraic point of view. So let $B$ be a regular local ring and $G$ a finite cyclic group of order $n$ acting faithfully on $B$ by local automorphisms. In the tame case; i.e. the order of $G$ is prime to the characteristic of the residue field $k$ of $B$, there is a central result of J.P. Serre \cite {serre} saying that the action is given by multiplying a suitable system of parameters $(y_1,\ldots,y_d)$ by roots of unity $y_i\longmapsto\zeta^{n_i}\cdot y_i$ for $i=1,\ldots, d$ where $\zeta$ is a primitive $n^{th}\!$-root of unity. Moreover, the ring of invariants $A:=B^G$ is regular if and only if $n_i\equiv 0 \mod n$ for $d-1$ of the parameters. The latter is equivalent to the fact that $\rk((\sigma-\id)|T)\le 1$ for the action of $\sigma\in G$ on the tangent space $T:={\frak m}_B/{\frak m}_B^2$.
 For more details see \cite[Chap. 5, ex. 7]{B-Lie}.\\
Only very little is known in the case of a wild group action; i.e., $\gcd(n,{\rm char}\;k)>1$. In this paper we will restrict ourselves to the case of $p\!$-cyclic group actions; i.e. $n=p$ is a prime number. We will present a sufficient condition for the fact that the ring of invariants $A$ is regular.  Our result is also valid in the tame case; i.e. where $n$ is a prime different from ${\rm char}\;k$. As the method of Serre depends on an intrinsic formula for writing down the action explicitly, we provide also an explicit formula for presenting $B$ as a free $A\!$-module if our condition is fulfilled. \\
The interest in our problem stems from the investigation on the relationship between the regular and the stable $R\!$-model of a smooth projective curve $X_K$ over the field of fractions $K$ of a discrete valuation ring $R$. In general, the curve $X_K$ admits a stable model $X'$ over a finite Galois extension $R\to R'$. Then the Galois group $G=G(R'/R)$ acts on $X'$. Our result provides a means to construct a regular model over $R$ starting from the stable model $X'$. We intend to work this out in a further article.\\[2mm]
{\small Finally, we want to mention that S. Wewers obtained partial results of ours by different methods cf. \cite{W}}\\

In this paper we will study only local actions of a cyclic group $G$ of prime order $p$ on a normal local ring $B$.
We fix a generator $\sigma$ of $G$ and obtain the {\it augmentation map}
$$
I:=I_{\sigma}:=\sigma-\id:B\To B\;;\;b\longmapsto \sigma(b)-b\;\;.
$$
We introduce the $B\!$-ideal
$$
I_{G}:= (I(b)\;;\;b\in B)\subset B
$$
which is generated by the image $I(B)$. This ideal is called {\it augmentation ideal}. If this ideal is generated by an element $I(y)$, we call $y$ an {\it augmentation generator.} Note that this ideal does not depend on the chosen generator $\sigma$ of $G$. Moreover, if $y$ is an augmentation generator with respect to a generator $\sigma$ of $G$, then $y$ is also an augmentation generator for any other generator of $G$. Since $B$ is local, the ideal $I_G$ is generated by an augmentation generator if $I_G$ is principal. Namely, $I_G/{\frak m}_BI_G$ is a vector space over the residue field $k_B=B/{\frak m}_B$ of $B$ of dimension $1$. So it is generated by the residue class of $I(y)$ for some $y\in B$ and, hence, due to Nakayama's Lemma, $I_G$ is generated by $I(y)$.

\begin{Def}\em
An action of a group $G$ on a regular local ring $B$ by local automorphisms is called a {\it pseudo-reflection} if there exists a system of parameters $(y_1,\ldots,y_d)$ of $B$ such that $y_2,\ldots,y_d$ are invariant under $G$.
\end{Def}

\begin{Thm}\label{Thm:MainThm}
 Let $B$ be a normal local ring with residue field $k_B:=B/{\frak m}_B$. Let $p$ be a prime number and
$G$ a $p\!$-cyclic group of local automorphisms of $B$.
Let $I_{G}$ be the augmentation ideal.
Let $A$ be the ring of $G\!$-invariants of $B$.
Consider the following conditions:\\[2mm]
\itboxx{a} $I_G:=B\cdot I(B)$ is principal.\\[2mm]
\itboxx{b} $B$ is a monogenous $A\!$-algebra.\\[2mm]
\itboxx{c} $B$ is a free $A\!$-module.\\[2mm]
Then the following implications are true:
$$
({\rm a})\longleftrightarrow ({\rm b})\longrightarrow ({\rm c})
$$
Assume, in addition, that $B$ is regular. Consider the following conditions:\\[2mm]
\itboxx{d} $A$ is regular.\\[2mm]
\itboxx{e} $G$ acts as a pseudo-reflection.\\[2mm]
Then the  condition {\rm (c)} implies {\rm (d)}.\\[2mm]
Moreover if, in addition, the canonical map $k_A\,\tilde{\To}\,k_B$ is an isomorphism. Then
condition  {\rm (a)} is equivalent to condition {\rm (e)}.
\end{Thm}

We start the proof of the theorem by several preparations.

\begin{Rem}\label{Rem:I-map}
\em For $b_1,b_2,b\in B$, the following relations are true:\\[2mm]
\itboxx{i} $I(b_1\cdot b_2)=I(b_1)\cdot\sigma(b_2)+b_1\cdot I(b_2)$\\[1mm]
\itboxx{ii} $I(b^n)=\ds \left(\sum_{i=1}^{n}\sigma(b)^{i-1}b^{n-i}\right)\cdot I(b)$\\[1mm]
\itboxx{iii} $\ds I\left( \frac{b_1}{b_2} \right)=\frac{I(b_1)b_2-b_1I(b_2)}{b_2\sigma(b_2)}$ if $b_2\ne 0$.
\end{Rem}
{\it Proof.} (i) follows by a direct calculation and (ii) by induction from (i).\\
(iii) The formula (i) holds for elements in the field of fractions as well. Therefore it holds
$$
I(b_1)=I\left( \frac{b_1}{b_2}b_2 \right)=I\left( \frac{b_1}{b_2} \right)\sigma(b_2)+\frac{b_1}{b_2}I(b_2)
$$
and the formula follows.\ok

For the implication (a)$\to$(b) we need a technical lemma.

\begin{Lemma} \label{Lemma:Main}
Let $y\in B$ be an augmentation generator. Then set, inductively,
$$
\begin{array}{lcll}
y_{i}^{(0)}&:=&y^{i}&\;\;\mbox{for}\;\;i=0,\ldots,p-1\\[2mm]
y_{i}^{(1)}&:=&\ds\frac{I\left(y_{i}^{(0)}\right)}{I\left(y_{1}^{(0)}\right)} & \;\;\mbox{for}\;\;i=1,\ldots,p-1\\[5mm]
y_{i}^{(n+1)}&:=&\ds\frac{I\left(y_{i}^{(n)}\right)}{I\left(y_{n+1}^{(n)}\right)}& \;\;\mbox{for}\;\;i=n+1,\ldots,p-1\;\;.
\end{array}
$$
Then
$$
y_{i}^{(n)}=\sum_{0\le k_1\le\ldots\le k_{i-n}\le n}\;\prod_{j=1}^{i-n}\sigma^{k_j}(y)\;\;\mbox{for}\;\;i=n,\ldots,p-1
$$
and, in particular,
$$
\begin{array}{lcl}
y_{n}^{(n)}&=&\ds1\\[2mm]
y_{n+1}^{(n)}&=&\ds \sum_{j=1}^{n+1}\sigma^{j-1}(y) \\[2mm]
I\left( y_{n+1}^{(n)} \right)&=&\sigma^{n+1}(y)-y\\[2mm]
\end{array}
$$
Furthermore, $y_{n+1}^{(n)}$ is again an augmentation generator for $n=0,\ldots,p-2$.
\end{Lemma}
{\it Proof.} We proceed by induction on $n$.
For $n=0$ the formulae are obviously correct.
For the convenience of the reader we also display the formulae for  $n=1$.  Due to Remark~\ref{Rem:I-map} one has
$$
\begin{array}{lcl}
y_i^{(1)}&=&\ds\frac{I\left( y_{i}^{(0)} \right)}{I\left(y_1^{(0)}\right)}=\frac{I(y^{i})}{I(y)}=\sum_{j=1}^{i}\sigma(y)^{j-1}y^{i-j}\\[2mm]
&=&\ds\sum_{0\le k_1\le\ldots\le k_{i-1}\le 1}\;\prod_{\nu=1}^{i-1}\sigma^{k_\nu}(y)
\end{array}
$$
since the last sum can be viewed as a sum over an index $j$ where $i-j$ is the number of the $k_{\nu}=0$. In particular, the formulae are correct for $y^{(1)}_1$ and $y_{2}^{(1)}$. Moreover
$$
I\left( y_{2}^{(1)} \right)=I(\sigma(y)-y)=\sigma^{2}(y)-y\;\;.
$$
Since $\sigma^{2}$ is generator of $G$ for $2<p$, the element $y_{2}^{(1)}$ is an augmentation generator as well.\\
Now assume that the formulae are correct for $n$. Since $y_{n+1}^{(n)}$ is an augmentation generator, $I\left( y_{n+1}^{(n)} \right)$ divides $I\left(  y_{i}^{(n)} \right) $ for $i=n+1,\ldots ,p-1$. Then it remains to show
$$
I\left(  y_{i}^{(n)} \right)=\left( \sigma^{n+1}(y)-y \right)\cdot  y_{i}^{(n+1)}\;\;\mbox{for}\;\;i=n+1,\ldots,p-1\;\;.
$$
For the left hand side one computes
$$
\begin{array}{lcl}
LHS&=&\ds I\left(
\sum_{0\le k_1\le\ldots\le k_{i-n}\le n}\;\prod_{j=1}^{i-n}\sigma^{k_j}(y)
 \right) =
 \sum_{0\le k_1\le\ldots\le k_{i-n}\le n}\;I
 \left(
  {\prod_{j=1}^{i-n}\sigma^{k_j}(y)}
  \right)\\[5mm]
  &=&\ds \sum_{0\le k_1\le\ldots\le k_{i-n}\le n}\cdot
 \left(
 \prod_{j=1}^{i-n}\sigma^{k_j+1}(y)-\prod_{j=1}^{i-n}\sigma^{k_j}(y)
  \right)\\[5mm]
  &=&
  \ds \sum_{1\le k_1\le\ldots\le k_{i-n}\le n+1}\;\prod_{j=1}^{i-n}\sigma^{k_j}(y)
 \; - \sum_{0\le k_1\le\ldots\le k_{i-n}\le n}\;\prod_{j=1}^{i-n}\sigma^{k_j}(y)\;\;.
\end{array}
$$
Now all terms occurring in both sums cancel. These are the terms with $k_{i-n}\le n$ in the first sum and $1\le k_1$  in the second sum.\\
For the right hand side one computes
$$
\begin{array}{lcl}
RHS&=&\ds \left( \sigma^{n+1}(y)-y \right)\cdot\sum_{0\le k_1\le\ldots\le k_{i-n-1}\le n+1}\;\prod_{j=1}^{i-n-1}\sigma^{k_j}(y)\\[5mm]
&=&\ds\sum_{0\le k_1\le\ldots\le k_{i-n}= n+1}\;\prod_{j=1}^{i-n}\sigma^{k_j}(y)
\;-\sum_{0= k_1\le\ldots\le k_{i-n}\le n+1}\;\prod_{j=1}^{i-n}\sigma^{k_j}(y)\;\;.
\end{array}
$$
Comparing both sides one obtains $LHS=RHS$. In particular we have
$$
\begin{array}{lcl}
y_{n+1}^{(n+1)}&=&1\\[2mm]
y_{n+2}^{(n+1)}&=&\ds\sum_{0\le k_1\le n+1}\prod_{j=1}^{1}\sigma^{k_1}(y)=\sum_{j=1}^{n+2}\sigma^{j-1}(y)\\[5mm]
I\left( y_{n+2}^{(n+1)} \right)&=&\sigma^{n+2}(y)-y\;\;.
\end{array}
$$
So $y_{n+2}^{(n+1)}$ is an augmentation generator for $n+2<p$, since $\sigma^{n+2}$ generates $G$.
This concludes the technical part.\ok

\begin{Prop} \label{Prop:MainProp}
Assume that the augmentation ideal $I_G$ is principal and let $y\in B$ be an augmentation generator. Then $B$ decomposes into the direct sum
$$
B=A\cdot y^0\oplus A\cdot y^1\oplus\ldots\oplus A\cdot y^{p-1}\;\;.
$$
\end{Prop}
{\it Proof.} Since $I(y)\ne 0$, the element $y$ generates the field of fractions $Q(B)$ over $Q(A)$.
Therefore
$$
Q(B)=Q(A)\cdot y^0\oplus Q(A)\cdot y^1\oplus\ldots\oplus Q(A)\cdot y^{p-1}\;\;.
$$
Then it suffices to show the following claim:\\[2mm]
Let $a,a_0,\ldots,a_{p-1}\in A$. Assume that $a$ divides
$$
b=a_0\cdot y^0+a_1\cdot y^1+\ldots +a_{p-1}\cdot y^{p-1}\;\;.
$$
Then $a$ divides $a_0,a_1,\ldots,a_{p-1}$.\\[2mm]
If $b=a\cdot \beta$, then $I(b)=a\cdot I(\beta)$. Since $I(\beta)=\beta_1\cdot I(y)$, we get $I(b)=a\beta_1\cdot I(y)$. So we see that $a$ divides $I(b)/I(y)\in B$. Using the notations of Lemma~\ref{Lemma:Main}, set
$$
\begin{array}{lclcl}
b^{(0)}&:=&b&=&a_0\cdot y^0+a_1\cdot y^1+\ldots a_{p-1}\cdot y^{p-1}\\[2mm]
b^{(1)}&:=&\ds\frac{I\left(b^{(0)}\right)}{I(y)}&=&\ds a_1+a_2\frac{I\left(y^{2}\right)}{I(y)}+\ldots+a_{p-1}\frac{I\left(y^{p-1}\right)}{I(y)}\\[4mm]
&&&=&a_1\cdot y_1^{(1)}+a_{2}\cdot y_{2}^{(1)}+\ldots +a_{p-1}\cdot y_{p-1}^{(1)}\\[2mm]
b^{(n)}&:=&\ds\frac{I\left(b^{(n-1)}\right)}{I\left(y^{(n-1)}_{n}\right)}&=&a_n\cdot y_n^{(n)}+a_{n+1}\cdot y_{n+1}^{(n)}+\ldots +a_{p-1}\cdot y_{p-1}^{(n)}\;\;.
\end{array}
$$
Due to the observation above, we see by induction that $a$ divides $b^{(0)},b^{(1)},\ldots,b^{(p-1)}$, since $y^{(n)}_{n+1}$ is an augmentation generator for $n=1,\ldots, p-2$. So we obtain
$$
a\,\big|\,b^{(p-1)}=a_{p-1}\cdot y^{(p-1)}_{p-1}=a_{p-1}\;\;.
$$
Now proceeding downwards, one obtains
$$
\begin{array}{l}
a\,\big|\,b^{(p-2)}=a_{p-2}+a_{p-1}\cdot y^{(p-2)}_{p-1}\;\;\mbox{and, hence,}\;\;a\,\big|\,a_{p-2}\\[2mm]
a\,\big|\,b^{(n)}=a_{n}+a_{n+1}\cdot y^{(n)}_{n+1}+\ldots+a_{p-1}\cdot y^{(n)}_{p-1}\;\;\mbox{and, hence,}\;\;a\,\big|\,a_{n}
\end{array}
$$
for $n=p-1,p-2,\ldots,0$.\ok

{\it Proof of the first part of Theorem}~\ref{Thm:MainThm}.\\[2mm]
(a)$\to$(b): This follows from Proposition~\ref{Prop:MainProp}.\\
(b)$\to$(a): If $B=A[y]$ is monogenous, then $I_G=B\cdot I(y)$ is principal.\\
(b)$\to$(c) is clear. Namely, if $B=A[y]$, the minimal polynomial of $y$ over the field of fraction is of degree $p$ and the coefficients of this polynomial belong to $A$. Then $B$ has $y^0,y^1,,\ldots,y^{p-1}$ as an $A\!$-basis. \\

Next we do some preparations for proving the second part of the theorem where $B$ is assumed to be regular.

\begin{Lemma}\label{Lemma:I-gen}
Let $R$ be a local subring of $B$ which is invariant under $G$ such that the canonical map $R/{\frak m}_R\,\tilde{\,\To}\,B/{\frak m}_B$ is an isomorphism. Let $(y_1,\ldots,y_d)$ be a generating system of the maximal ideal ${\frak m}_B$. Then $I_G$ is generated by $I(y_1),\ldots,I(y_d)$.
\end{Lemma}
{\it Proof.} Due to the assumption, we have $B=R+{\frak m}_B$ and, hence, $I(B)=I({\frak m}_B)$. Furthermore, we have
$$
{\frak m}_B={\frak m}_B^2+\sum_{i=1}^{d}R\cdot y_i\;\;.
$$
Since $I$ is $R\!$-linear, we get
$$
I({\frak m}_B)=I\left({\frak m}_B^2\right)+\sum_{i=1}^{d}R\cdot I(y_i)\;\;.
$$
Due to Remark~\ref{Rem:I-map}, one knows $I\left({\frak m}_B^2\right)\subset {\frak m}_B\cdot I({\frak m}_B)$. So one obtains
$$
I({\frak m}_B)\subset {\frak m}_B\cdot I({\frak m}_B)+\sum_{i=1}^{d}R\cdot I(y_i)\;\;.
$$
Since $B$ is local, Nakayama's Lemma yields
$$
I_G=B\cdot I(B)=B\cdot I({\frak m}_B)=\sum_{i=1}^{d}B\cdot I(y_i)\;\;.
$$
Thus the assertion is proved.\ok

\begin{Prop} \label{Prop:I-gen}
Keep the assumption of the second part of the theorem; namely that $B$ is regular and that the canonical morphism  $k_A\,\tilde{\To}\,k_B$ is an isomorphism. Let $(y_1,\ldots,y_d)$ be a generating system of the maximal ideal ${\frak m}_B$. Then the following assertions are true:\\[2mm]
\itboxx{i} $I_G=B\cdot I(y_1)+\ldots+B\cdot I(y_d)$ \\[2mm]
\itboxx{ii} If the ideal $I_G=B\cdot I(B)$ is principal, then there exists an index $i\in\{1,\ldots,d\}$ with $I_G=B\cdot I(y_{i})$.
\end{Prop}
{\it Proof.} Let $\widehat{B}$ be the ${\frak m}_B\!$-adic completion of $B$. Recall that  a regular ring is noetherian by definition.
Therefore the extension $B\to \widehat{B}$ is faithfully flat and ${\frak m}_B\widehat{B}={\frak m}_{\widehat{B}}$ ; cf. \cite[10.14 \& 10.15]{AM}. Since $G$ acts by local morphism, any $\sigma\in G$ extends to a local automorphism $\widehat{\sigma}$ of $\widehat{B}$.
Due to the assumption that the canonical morphism  $k_A\,\tilde{\To}\,k_B$ is an isomorphism, any $b\in B$ can be written as $B=a+m$ where $a\in A$ is invariant under $G$ and $m\in {\frak m}_B$ and, hence, $I(b)=I(m)\in I({\frak m}_B)$. Therefore $I_G$ is generated by $I({\frak m}_B)$.\\
(i) Since $B\to \widehat{B}$ is faithfully flat and $\widehat{B}\cdot {\frak m}_B=\frak m_{\widehat{B}}$, it suffices to prove the assertion for the completion $\widehat{B}$. For complete local rings there exists a $G\!$-stable lift $R$ of the residue field $k$. Namely, in the case of mixed characteristic $(0,p)$, one can choose the ring of Witt vectors $W(k)\subset \widehat{A}$ as $R$ and, in the equal characteristic case, the residue field $k$ lifts into $\widehat{A}$; cf. \cite{C}. Now we can apply Lemma~\ref{Lemma:I-gen} and obtain the assertion.\\
(ii) Since $I_G$ is principal, $I_G/{\frak m}_BI_G$ is generated by one of the $I(y_i)$ and, hence, again by Nakayama's Lemma
$I_G=B\cdot I(y_i)$ for a suitable $i\in\{1,\ldots,d\}$.\ok

{\it Proof of the second part of Theorem}~\ref{Thm:MainThm}.\\[2mm]
(c)$\to$(d) follows from \cite[Theorem 51]{matsumura}. Namely, $B$ is noetherian due to the definition of a regular ring. Since $A\to B$ is faithfully flat, so $A$ is noetherian. Then one can apply loc.cit.\\
(a)$\to$(e) We assume that the canonical map $k_A\to k_B$ of the residue fields is an isomorphism.
If $I_G$ is principal, one can choose an augmentation generator $y\in {\frak m}_B$ which is part of a system of parameters $(y,y_2,\dots,y_d)$ due to Proposition~\ref{Prop:I-gen}. Due to Proposition~\ref{Prop:MainProp}, we know that $B$ decomposes into the direct sum
$$
B=A\cdot y^0\oplus A\cdot y^1\oplus\ldots\oplus A\cdot y^{p-1}\;\;.
$$
Now we can represent
$$
y_j=\sum_{i=0}^{p-1}a_{i,j}\cdot y^i\;\;\mbox{for}\;\;j=2,\ldots,d\;\;.
$$
Then set
$$
\tilde{y}_j:=y_j-\sum_{i=1}^{p-1}a_{i,j}y^i=a_{0,j}\in A\cap {\frak m}_B={\frak m}_A\;\;\mbox{for}\;\;j=2,\ldots,d\;\;.
$$
So $(y,\tilde{y}_2,\ldots,\tilde{y}_d)$ is a system of parameters of $B$ as well. Thus  $G$ acts by a pseudo-reflection.\\
(e)$\to$(a): If $G$ is a pseudo-reflection, $I_G$ is generated by $I(y)$ due to Proposition~\ref{Prop:I-gen} where $y,x_2,\ldots,x_{p}$ is a system of parameters with $x_i\in {\frak m}_A$ for $i=2,\ldots,p$ if $k_A = k_B$.\ok

If $k_A\to k_B$ is not an isomorphism, the implication (e)$\to$(a) is false as the following shows.

\begin{Ex}\em
Let $k$ be a field of positive characteristic $p$ and look at the polynomial ring
$$R:=k[Z,Y,X_1,X_2]$$
over $k$. We define a $p$-cyclic action of $G=\langle \sigma\rangle $ on $R$
by
$$\sigma|k:=\id_k\,,\,\sigma(Z)=Z+X_1\,,\, \sigma(Y)=Y+X_2\,,\,  \sigma(X_i)=X_i\;\;\mbox{for}\;\;i=1,2\;\;.$$
This is a well-defined action of order $p$, since $p\cdot X_i=0$ for $i=1,2$, and it leaves the ideal ${\frak I}:=(Y,X_1,X_2)$ invariant.
Furthermore, for any $g\in k[Z]-\{0\}$ the image is given by $\sigma(g)=g+I(g)$ with $I(g)\in X_1\cdot k[Z,X_1]$.\\
Then consider the polynomial ring
$$
S:=k(Z)[Y,X_1,X_2]
$$
over the field of fractions $k(Z)$ of the polynomial ring $k[Z]$. Then $S$ has the maximal ideal ${\frak m}=(Y,X_1,X_2)$. Then set
$$
B:=S_{{\frak m}}=k(Z)[Y,X_1,X_2]_{(Y,X_1,X_2)}\;\;.
$$
We can regard all these rings as subrings of the field of fractions of $R$
$$
R\subset S \subset B\subset k(Z,Y,X_1,X_2)\;\;.
$$
Clearly, $\sigma$ acts on $R$ and, hence, it induces an action on its field of fractions; denote this action by $\sigma$ as well.
Then we claim that the restriction of $\sigma$ to $B$ induces an action on $B$ by local automorphisms. For this, it suffices to show that
for any $g\in R-{\frak I}$ the image $\sigma(g)$ does not belong to ${\frak I}$. The latter is true, since
$$
\sigma(g)=g+I(g)\;\;\mbox{with}\;\;I(g)\in {\frak I}\;\;.
$$
The augmentation ideal $I_G=B\cdot X_1+B\cdot X_2$ is not principal although $G$ acts through a pseudo-reflection.\ok
\end{Ex}

\begin{Rem}\em
In the tame case $p\ne{\rm char}(k_B)$, the converse $({\rm d})\to ({\rm a})$ is also true due to the theorem of Serre as explained in the introduction.
\end{Rem}

In the case of a wild group action; i.e. $p={\rm char}(k_B)$, it is not known whether the converse is true, but we would conjecture it.

\begin{Conj}\em Let $B$ be a regular local ring and let $G$ be a $p\!$-cyclic group acting on $B$ by local automorphisms. Then the following conditions are {\it conjectured} to be equivalent:\\[2mm]
\itboxx{1} $I_G$ is principal.\\[2mm]
\itboxx{2} $A:=B^G$ is regular.
\end{Conj}
The implication $(1)\to (2)$ was shown in Theorem \ref{Thm:MainThm}. Of course the converse is true if $\dim A\le1$. In higher dimension, the converse $(2)\to (1)$ is uncertain, but it holds for small primes $p\le 3$ as we explain now.
Since $A$ is regular, the ring $B$ is a free $A\!$-module of rank $p$; cf. \cite[IV, Prop. 22]{serre2}. So,
$$
B/B{\frak m}_A^n\;\;\mbox{is a free $A/{\frak m}_A^n\!$-module of rank}\;\; p\;\;\mbox{for any}\;\;n\in \mathbb{N}\;\;.
\leqno(*)
$$
In the case $p=2$ the rank of ${\frak m}_B/B{\frak m}_A$ is $0$ or $1$. In the first case, $k_B$ is an extension of degree $[k_B:k_A]=2$ over $k_A$ and ${\frak m}_B=B{\frak m}_A$. So there exists an element $\beta\in B$ such that $B/B{\frak m}_A$ is generated by the residue classes of $1$ and $\beta$. Due to Nakayama's Lemma $B=A[\beta]$ is monogenous and, hence, $I_G$ is principal. In the second case, where $k_A\to k_B$ is an isomorphism, then there exists an element $\beta\in {\frak m}_B$ such that ${\frak m}_B=B\beta+B{\frak m}_A$. Then $G$ acts as a pseudo-reflection and, hence, $I_G$ is principal.\\[2mm]
In the case $p=3$ we claim that $B{\frak m}_A\not\subset{\frak m}_B^2$!\\
If we assume the contrary  $B{\frak m}_A\subset{\frak m}_B^2$ then these ideals coincide; $B{\frak m}_A={\frak m}_B^2$. Namely, the rank of $B/B{\frak m}_A$ as $A/{\frak m}_A\!$-module is $3$ and the rank of $B/{\frak m}_B^2$ is at least $3$ due to $d:=\dim B\ge 2$, so $B{\frak m}_A={\frak m}_B^2$. Therefore the length of $B/B{\frak m}_A^2=B/{\frak m}_B^4$ is $3$ times the length of $A/{\frak m}_A^2$ which is $3\cdot(\dim A+1)$. On the other hand the rank of $B/{\frak m}_B^4$ is equal to
$$
\big(1+\dim {\frak m}_B/{\frak m}_B^2\big)+\dim {\frak m}_B^2/{\frak m}_B^3+\dim {\frak m}_B^3/{\frak m}_B^4=\sum_{n=0}^{3}{d+n-1\choose d-1}
$$
which larger than
$$
\big(1+\dim {\frak m}_A/{\frak m}_A^2\big)+ \big(1+\dim {\frak m}_A/{\frak m}_A^2\big) + \big(1+\dim {\frak m}_A/{\frak m}_A^2\big)
$$
since for $d\ge 2$ holds
$$
{d+1\choose d-1}=\frac{(d+1)d}{2}\ge1+d=1+\dim {\frak m}_A/{\frak m}_A^2
$$
and
$$
{d+3-1\choose d-1}=\frac{(d+2)(d+1)d}{2\cdot 3}> 1+ d
$$
Here we used the formula for the number $\lambda_{n,d}$ of monomials $T_1^{m_1}\ldots T_d^{m_d}$in $d$ variables of degree $n=m_1+\ldots+m_d$
$$
\lambda_{n,d}={d+n-1 \choose d-1}\;\;.
$$
So, using only the condition $(*)$ and proceeding by induction on $\dim(A)$, we see that here exists a system of parameters $\alpha_1,\ldots,\alpha_d$ of $A$ such that $\alpha_2,\ldots,\alpha_d$ is part of a system of parameters of $B$. In the case, where $k_A\to k_B$ is an isomorphism, $G$ acts as a pseudo-reflection and, hence, $I_G$ is principal. If $k_A\to k_B$ is not an isomorphism, then we must have ${\frak m}_B=B{\frak m}_A$; otherwise the rank of $B/{\frak m}_B$ is at least $4$. Since $[k_B:k_A]\le 3$, the field extension $k_A\to k_B$ is monogenous and, hence, $A\to B$ is monogenous due to the Lemma of Nakayama.\ok

%\newpage
%\input{CGA-Conject}
{

} 

\begin{thebibliography}{BGR}
\itemsep0exminus0.1ex
\bibitem[AM]{AM}Atiyah, M.F.; Macdonald, I.G.: {\sl Commutative Algebra\/}. Addison Wesley Publishing Company, 1969.
\bibitem[B]{B-Lie}Bourbaki, N.: {\sl Groupes et alg\`{e}bre de Lie\/}. Masson, Paris, 1981
\bibitem[BGR]{BGR}Bosch, S.; Günzer, U.; Remmert, R.: {\sl Non-Archimedean Analysis\/}. Grundlehren {\bf 261}, Springer-Verlag, 1984.  
\bibitem[C]{C} Cohen, I.S.: {\sl On the structure and ideal theory of complete local rings\/}. Trans. Amer.
Math. Soc. {\bf 59}, 54-106 (1946).
\bibitem[E]{E}Edixhoven, B.: {\sl N\`{e}ron models and tame ramification\/}. Compsitio Math. {\bf 31}, 291-306, 1992.
\bibitem[H]{Hirzebruch}Hirzebruch, F.: {\sl Über vierdimensionale Riemannsche Flächen mehrdeutiger analytischer Funktionen von zwei komplexen Veränderlichen\/}. Math. Ann. {\bf 126}, 1-22, 1953.
\bibitem[J]{Jung} Jung, H. W. E.: {\sl Darstellung der Funktionen eines algebraischen Körpers zweier unabhängiger Veränderlicher $x,y$ in der Umgebung einer Stelle $x=a,y=b$\/}. J. Reine Angew. Math. {\bf 133}, 289-314, 1908.
\bibitem[M]{matsumura}Matsumura, H.: {\sl Commutative Algebra\/}. Benjamin, Mathematics Lecture Notes Series {\bf 56}, 1980.
\bibitem[S1]{serre}Serre, J.-P.: {\sl Groupes finis d'automorphismes d'anneaux locaux r\'{e}guliers\/}. In Colloque d'Alg\'{e}bre (Paris, 1967), Exp. 8. Secr\'{e}tariat math\'{e}matique, 1967.
\bibitem[S2]{serre2}Serre, J.-P.: {\sl Alg\`{e}bre Locale - Multiplicit\'{e}s\/}. Lecture Notes in Math. {\bf 11}, Springer-Verlag, 1965.
\bibitem[W]{W} Wewers, S.: {\sl Regularity of quotients by an automorphism of order $p$\/}. ArXiv:1001.0607, 2010.
\end{thebibliography}
\end{document}